\documentclass[3p,twocolumn]{elsarticle}               

\usepackage{graphicx}          
\usepackage{amsmath}
\usepackage{amssymb}
\usepackage{colonequals}

\newtheorem{definition}{Definition}
\newtheorem{thm}{Theorem}
\newtheorem{prop}[thm]{Proposition}
\newtheorem{rem}[definition]{Remark}
\newtheorem{lem}[thm]{Lemma}
\newtheorem{cor}[thm]{Corollary}
\newcommand{\sdcHide}[1]{{\bf #1}}
\renewcommand{\sdcHide}[1]{}

\journal{Automatica}

\begin{document}
	
	\begin{frontmatter}
		
		\title{MPC on manifolds with an application to the control of spacecraft attitude on $\textrm{SO}(3)$} 

		\fntext[footnoteinfo]{Preliminary results related to this research 
			were presented in~\cite{kalabicgupta_acc14}. Corresponding 
			author 
			U.~V.~Kalabi\'{c}. Tel. +1-617-621-7548. 
			Fax +1-617-621-7550.}
		\fntext[melco]{This work was not supported by Mitsubishi Electric Co.~or any of its subsidiaries.}

		\author[AnnArborAero]{Uro\v{s} V.~Kalabi\'{c}\fnref{footnoteinfo}}\ead{kalabic@umich.edu}    
		\author[AnnArborAero]{Rohit Gupta}\ead{rohitgpt@umich.edu}              
		\author[MERL]{Stefano Di Cairano\fnref{melco}}\ead{dicairano@ieee.org}
		\author[AnnArborMath]{Anthony M.~Bloch}\ead{abloch@umich.edu}
		\author[AnnArborAero]{Ilya V.~Kolmanovsky}\ead{ilya@umich.edu}

		\address[AnnArborAero]{Department of Aerospace Engineering, University of Michigan, Ann Arbor, MI 48109}  
		\address[MERL]{Mechatronics, Mitsubishi Electric Research Laboratories, Cambridge, MA 02139}     
		\address[AnnArborMath]{Department of Mathematics, University of Michigan, Ann Arbor, MI 48109}  
		
		

		\begin{abstract}                          
			We develop a model predictive control (MPC) design for 
			systems with discrete-time dynamics evolving on smooth 
			manifolds. 
			We show that the properties of conventional MPC for dynamics 
			evolving on $\mathbb R^n$ are preserved and we develop a 
			design procedure for achieving similar properties. 
			We also demonstrate that for discrete-time dynamics  on 
			manifolds with Euler characteristic 
			not equal to 1, there do not exist globally stabilizing, 
			continuous control laws. 
			The MPC law is able to achieve global asymptotic stability on these manifolds, because the MPC law may be discontinuous.
			We apply the method to spacecraft attitude control, where the 
			spacecraft attitude evolves on the Lie 
			group $\text{SO}(3)$ and for which a continuous globally stabilizing 
			control law does not exist. In this case, the MPC law is
			discontinuous and achieves global stability.
		\end{abstract}
		
	\end{frontmatter}
	
	\section{Introduction}
	Conventional model predictive control (MPC) \cite{rawlings_mpc_book} is developed for and
	usually applied to systems whose discrete-time dynamics evolve 
	on the ``flat'' normed vector space $\mathbb R^n$. However, 
	the configuration spaces of some systems are smooth manifolds that are not diffeomorphic to $\mathbb R^n$. To design the prediction dynamics for such systems, finite-dimensional manifolds may be embedded in $\mathbb R^n$, and standard integration schemes employed to derive the discrete-time update equation, enforcing their evolution on the manifold through the use of equality constraints. However, standard integration 
	techniques do not preserve 
	symmetries for systems evolving on 
	manifolds, 
	and this results in the integration not correctly representing the 
	actual dynamics. Because of this, specific methods 
	for integrating system dynamics that evolve on manifolds 
	have been  developed in, for example, 
	\cite{munthekass99,marsdenwest01,hairer_numerical_book}.
	
	As a motivating example, consider the attitude dynamics of a 
	rigid body evolving on the manifold $\textrm{SO}(3)$. This manifold is a $3$-dimensional 
	space which can be embedded in $\mathbb R^n$ for $n \geq 5$ but is not diffeomorphic to $\mathbb R^n$ for any $n$. 
	The group $\textrm{SO}(3)$ is a particular case of a compact Lie group; 
	for mechanical systems whose 
	configuration space is a compact Lie group, the Lie group 
	variational integrator (LGVI) \cite{lee_diss} has been developed to obtain
	discrete-time update equations that preserve the underlying 
	group structure. Unlike ordinary integration schemes, the LGVI also preserves the conserved quantities of motion \cite{leemcclamr_cca05} and is therefore a more realistic prediction model. 
	Also unlike ordinary integration schemes, the LGVI equations of motion are significantly 
	different from conventionally used, discrete-time dynamics. For instance,
	the LGVI for spacecraft attitude gives an equation of the form \cite{leeleok08}, 
	\begin{subequations}\label{equ:so3_dyn_intro}
		\begin{align}
		g_{k+1} &= g_kf_k, \\
		f_{k+1}J-J f_{k+1}^{\text T} &=  Jf_{k} -f_{k}^{\text T}J+ {h^2}u_k, 
		\label{equ:so3_fdyn_intro}
		\end{align}
	\end{subequations}
	where $g_k \in \textrm{SO}(3)$ is the 
	orientation, $f_k \in \textrm{SO}(3)$ is the change in $g_k$, $u_k$ is related to 
	the external torques, $h$ is the length of the integration step. Indeed, \eqref{equ:so3_fdyn_intro} is an implicit equation, and implicit equations are seldom used in conventional MPC prediction models.  
	
	The control of dynamics evolving on manifolds present 
	additional challenges. 
	On manifolds with Euler characteristic not equal to 1, there exist topological obstructions that imply the non-existence of globally-stabilizing, 
	continuous control laws. This has been shown for the continuous-time setting \cite{bhatbernstein00}, and is a consequence of the Poincar\'{e}-Hopf theorem \cite{guillemin_book}. In this paper, we appeal to the Lefschetz-Hopf theorem and  derive a similar result for discrete-time dynamics.
	

	In this paper we investigate 
	the design 
	of MPC for systems evolving in manifolds, with the aim of retaining  
	properties of conventional MPC applied to flat spaces. The applications of
	MPC to systems whose dynamics 
	evolve on manifolds have scarcely been considered, with the 
	exceptions of \cite{nagyfindeis_acc00,gros_nmpc12}, 
	which  focus on computational issues, and have limited 
	analysis of closed-loop properties. 
	
	In this paper\footnote{
		This work is a significant extension of the authors preliminary 
		investigation in~\cite{kalabicgupta_acc14}, which focused only 
		on  
		$\textrm{SO}(3)$, and did not present detailed proofs. The 
		results in this paper apply to
		general manifolds, and the theoretical 
		results are rigorously derived, including the non-existence of 
		discrete-time continuous 
		control laws on certain types of manifolds.}, we first show that 
	MPC on smooth manifolds
	achieves  the same properties of conventional
	MPC  in $\mathbb R^n$, namely, recursive feasibility and 
	asymptotic 
	stability of the equilibrium. Then we  describe the design of the 
	MPC
	terminal cost and terminal set that achieves such properties. 

	As an additional contribution of this paper, we show that when the manifold is compact, 
	our MPC is able to achieve global 
	closed-loop asymptotic stability in the state-unconstrained case.
	In particular, this property holds for manifolds with Euler characteristic not equal to 1, for which a 
	continuous, globally stabilizing control law does not exist.
	Because MPC can produce discontinuous control 
	laws \cite{meadows95}, the topological obstruction for such manifolds is not restrictive in this case.
	Thus this paper highlights the capability to design possibly discontinuous, stabilizing feedback control laws with global stability properties, 
	by employing a systematic and unified design procedure in the MPC  framework. 
	
	A practical contribution of this paper is the presentation of a 
	control law for the constrained 
	control of spacecraft attitude, in which we apply our 
	MPC scheme to the LGVI dynamics evolving on $\textrm{SO}(3)$. 

	The paper is organized as follows. Section \ref{sec:mpc} 
	develops the MPC scheme and the conditions for closed-loop 
	stability. 
	Section \ref{sec:locallaw} describes the design of the  terminal 
	penalty and terminal 
	set constraint achieving closed-loop stability
	by exploiting a specifically designed 
	local control law. 
	Section \ref{sec:globallaw} proves the non-existence of a globally 
	stabilizing, continuous control law under certain assumptions, 
	and hence the discontinuity of the MPC law. 
	Section \ref{sec:so3} presents simulation results for dynamics 
	evolving on the Lie group 
	$\textrm{SO}(3)$ and  Section \ref{sec:conc} summarizes the 
	conclusions.

	\sdcHide{(is it necessary to talk about geometric control here?) 
		The 
		field of geometric, nonlinear control is a mature field with various 
		topics that have been extensively addressed (\text{e.g.}, see 
		\cite{jurdjevic_book,bloch_book,agrachev_book,bullo_book} }

	\sdcHide{( this must go somewhere else --probably when you 
		describe 
		the results for DT): The fixed point index depends 
		on the system 
		dynamics, but is invariant a long as the dynamic update 
		equations are homotopic to the identity map. The fixed point 
		index corresponding to a sink is 1, so it follows that if the Euler 
		characteristic is not equal to 1, there exists no map such that the 
		sole equilibrium is a sink.}

	
	\subsection{Notation}
	The notation is standard with a few notable 
	exceptions. The set $\mathbb Z_N$ denotes the set of the first 
	$N$ nonnegative integers and $\mathbb Z_+$ denotes the set of 
	all nonnegative integers. For a set $\mathcal A$, its interior is 
	denoted by $\operatorname{int}\mathcal A$, its closure by 
	$\operatorname{cl}\mathcal A$, and its boundary by 
	$\operatorname{bd}\mathcal A$. The set $\mathcal A^N$ 
	represents $N$ copies of $\mathcal A$. The identity matrix is 
	denoted $ I_n \in \mathbb R^{n \times n}$ and an $n$-by-$m$ 
	zero matrix is denoted $0_{n \times m} \in \mathbb R^{n \times 
		m}$. For two functions $f, g: \mathbb R \rightarrow \mathbb R$, $f \circ g$ is their composition, and
	$g(t) = o(f(t))$ implies that $\lim_{t \rightarrow \infty} g(t)/f(t) = 
	0$. 
	A function $\alpha: [0,a) \rightarrow [0,\infty)$ is said to be class 
	$\mathcal K$ if it is strictly increasing and $\alpha(0) = 0$; 
	furthermore $\alpha$ is said to be class $\mathcal K_\infty$ if $a 
	= \infty$ and $\alpha(x) \rightarrow \infty$ as $x \rightarrow 
	\infty$.
	Finally, for a sequence $\{v_k,v_{k+1},\dots,v_{k+N}\}$, its predicted value at time $k$ is denoted by $\{v_{k|k},v_{k+1|k},\dots,v_{k+N|k}\}$.
	
	
	\section{MPC on manifolds}\label{sec:mpc}
	
	We begin by developing a general MPC law for application to 
	smooth manifolds. In the construction of the control law, we follow ideas inspired by conventional nonlinear MPC laws 
	\cite{rawlings_mpc_book}, 
	for systems 
	whose dynamics evolve on $\mathbb R^n$. For instance, in order to ensure 
	recursive feasibility of the finite horizon optimal control problem and to enlarge the domain 
	of attraction, we utilize a terminal 
	set and terminal cost function \cite{rawlings_mpc_book}. 
	We thereby obtain a result on system stability that rigorously generalizes MPC stability results to a more generous class of systems whose dynamics evolve on manifolds.
	In this setting, we obtain a remarkable result that when the manifold is compact and there are no 
	state constraints, this MPC law may provide global asymptotic stability.
	
	Let $M$ be an $n$-dimensional smooth manifold, which is endowed with a metric $d$, and let $U$ be a compact subset of an $m$-dimensional smooth manifold. Note that the fact that all smooth manifolds are metrizable follows from Whitney's embedding theorem \cite{guillemin_book}. Furthermore, all smooth manifolds admit a Riemannian metric \cite{lee_book}. Consider the dynamic update equation,
	\begin{subequations}\label{equ:dyn}
		\begin{gather}
		x_{k+1} = f(x_k,u_k), \\
		f: M \times U \rightarrow M, \label{equ:dyn2}
		\end{gather}
	\end{subequations}
	where $x_k \in M$ and $u_k \in U$. The evolution on manifold is highlighted by using $M$ in the domain and co-domain of \eqref{equ:dyn2}; the function $f$ is continuously differentiable and satisfies $f(x_e,u_e) = x_e$ for some $x_e \in M$, which we refer to as the equilibrium of $f$, and some $u_e \in U$.
	
	The system is subject to state and control constraints,
	\begin{equation}\label{equ:constr}
	x_k \in \mathcal X,~u_k \in \mathcal U,
	\end{equation}
	where $\mathcal X$ and $\mathcal U$ are compact and connected subsets of $M$ and $U$, respectively, that satisfy $x_e \in \operatorname{int}\mathcal X$ and $u_e \in \operatorname{int}\mathcal U$. 
	
	We introduce a cost function $\mathcal V_N: M \times U^N \rightarrow \mathbb R$ satisfying,
	\begin{multline}\label{equ:mpccost}
	\mathcal V_N\left(x_k;\{u_{k+i|k}\}_{i \in \mathbb Z_N}\right) \\ = F(x_{k+N|k})+\sum_{i \in \mathbb Z_N} L(x_{k+i|k},u_{k+i|k}),
	\end{multline}
	where $x_{k+i+1|k} = f(x_{k+i|k},u_{k+i|k})$ is the predicted state at time $k+i+1$ for $i \in \mathbb Z_N$. 
	The functions $L: M \times U \rightarrow \mathbb R$ and $F: M \rightarrow \mathbb R$ are twice continuously differentiable and have the following properties,
	\begin{subequations}\label{equ:costprops}
		\begin{align}
		L(x_e,u_e) = F(x_e) = 0, & \label{equ:costprops1} \\
		L(x_k,u_k) \geq L(x_k,u_e) \geq \gamma(d(x_k,x_e)), 
		\label{equ:costprops2} \\
		F(x_k) \geq \alpha(d(x_k,x_e)), 
		\label{equ:costprops3}
		\end{align}
	\end{subequations}
	for all $x_k \in M$ and $u_k \in U$, where $\gamma$ is a class $\mathcal{K}_\infty$, $\alpha$ is a class $\mathcal{K}$ function. 
	
	
	We now introduce the target set and terminal feedback law. These are required in order to enforce recursive feasibility, ensuring that if the state at the end of the prediction horizon is in the target set, it and the associated feedback control satisfy all constraints and are therefore feasible solutions to the finite horizon optimal control problem. Specifically, we introduce the target set $\mathcal X_T \subset \mathcal X$, which is compact and contains $x_e$ in its interior. We also introduce a control law $\kappa: M \rightarrow U$, which we refer to as the local control law satisfying,
	\begin{subequations}\label{equ:locallaw}
		\begin{align}
		\kappa(x_k) &\in \mathcal U, \label{equ:locallaw1} \\
		f(x_k,\kappa(x_k)) &\in \mathcal X_T, \label{equ:locallaw2} \\
		F(f(x_k,\kappa(x_k)))-F(x_k) &\leq -L(x_k,\kappa(x_k)), \label{equ:locallaw3}
		\end{align}
	\end{subequations}
	for all $x_k \in \mathcal X_T$. 
	
	The MPC control law is obtained through the solution to the following problem,
	\begin{subequations}\label{equ:mpcopt}
		\begin{align}
		\min_{\{u_{k+i|k}\}_{i \in \mathbb Z_N}}~& \mathcal V_N\left(x_k;\{u_{k+i|k}\}_{i \in \mathbb Z_N}\right), \\  
		\operatorname{subject~to} ~
		& x_{k+i+1|k} = f(x_{k+i|k},u_{k+i|k}), \label{equ:mpcopt_1} \\
		& x_{k+i|k} \in \mathcal X,  \\
		& u_{k+i|k} \in \mathcal U,~\forall i \in \mathbb Z_N, \label{equ:mpcopt_2} \\
		& x_{k+N|k} \in \mathcal X_T. \label{equ:mpcopt_3}
		\end{align}
	\end{subequations}
	When the solution to (\ref{equ:mpcopt}) exists, it is denoted by $\mathcal V_N^*(x_k)$ and the control sequence solving it is denoted by $\{u_{k+i|k}^*\}_{i \in \mathbb Z_N}$. 
	
	The input obtained from the model predictive control law at time $k$ is the first element in the sequence solving (\ref{equ:mpcopt}),
	\begin{equation}\label{equ:mpclaw}
	u_k = u_{k|k}^*.
	\end{equation}
	
	The MPC law defined above can be used
	on differentiable manifolds. We will now show that 
	the domain of attraction of the equilibrium $x_e$ of the closed-loop system defined by (\ref{equ:dyn}) and (\ref{equ:mpclaw}) coincides with the set of initial conditions that can be steered using open-loop control to the target set without violating the constraints. 
	Define,
	\begin{multline}\label{equ:Dn_defn}
	\mathcal D_N = \{x_k \in M: \exists \{u_{k+i|k}\}_{i \in \mathbb Z_N} \in \mathcal U \\ \textrm{ s.t.~(\ref{equ:mpcopt_1})-(\ref{equ:mpcopt_3}) are satisfied}\}.
	\end{multline}
	We state the main result below and provide the proof in the appendix.
	
	\begin{thm}\label{thm:mpc}
		Let $x_0 \in \mathcal D_N$. Then (i) the finite horizon optimal control problem (\ref{equ:mpcopt}) is feasible for all $k \in \mathbb Z_+$ and (ii) $x_e \in M$ is the asymptotically stable equilibrium for closed-loop system defined by (\ref{equ:dyn}) and (\ref{equ:mpclaw}) with a basin of attraction $\mathcal D_N$, \textit{i.e.}, for any $x_0 \in \mathcal D_N$ and $\varepsilon > 0$, there exists $\delta > 0$ such that $d(x_0,x_e) < \delta$ implies that $d(x_k,x_e) < \varepsilon$ for all $k > 0$, and furthermore, $d(x_k,x_e) \rightarrow 0$ as $k \rightarrow \infty$.
	\end{thm}
	
	Note that whenever $\mathcal D_N = M$, Theorem \ref{thm:mpc} implies that the control law (\ref{equ:mpclaw}) is globally stabilizing. In general, this may not be possible because $\mathcal D_N$ is a subset of $\mathcal X$ and $\mathcal X$ is a subset of $M$. 
	However, if $M$ is a compact manifold, then $\mathcal X$ may be equal to $M$; in this case, a large enough $N$ can be chosen so that $\mathcal D_N = \mathcal X$.
	Therefore, based 
	on the following proposition, it is possible to guarantee global asymptotic stability in the case where $M$ is compact. 
	
	
	\begin{prop}\label{prop:gas}
		Suppose $M$ is compact and let $\mathcal D_\infty$ be the set of initial conditions $x_0$ such that there exists a sequence of control inputs steering $x_0$ to $x_e$,
		\begin{multline}
		\mathcal D_\infty = \{x_k \in M: \exists \{u_{k+i|k}\}_{i \in \mathbb Z_+} \in \mathcal U, \\ \lim_{k \rightarrow \infty} x_{k+i|k} = x_e\}. 
		\end{multline}
		Suppose the system (\ref{equ:dyn}) is state-unconstrained so that $\mathcal X = M$ and suppose $\mathcal D_\infty = M$. Then there exists a finite $N^*$ such that $\mathcal D_{N^*} = M$ and therefore the closed-loop system defined by (\ref{equ:dyn}) and (\ref{equ:mpclaw}) is globally asymptotically stable.
	\end{prop}
	
	\textit{Proof.} It is clear that for any finite $N \geq 0$, $\mathcal D_N \subset \mathcal D_\infty = M$ and that the set $\operatorname{cl}(M \setminus \mathcal X_T) = \operatorname{cl}(\mathcal D_\infty \setminus \mathcal X_T)$ is a subset of $\mathcal D_\infty$. 
	Let $N(x) = \inf \{i: x_k = x,~\exists \{u_{k+i|k}\}_{i \in \mathbb Z_+} \in \mathcal U, ~x_{k+i|k} \in \mathcal X_T\}$ be the minimum number of control steps required to guide an initial condition $x \in D_\infty = M$ to the set $\mathcal X_T$. According to the definition of $\mathcal D_\infty$ and the fact that $x_e \in \operatorname{int}\mathcal X_T$, $N(x)$ is finite for all $x \in \operatorname{cl}(M \setminus \mathcal X_T)$. 
	Let $N^* = \sup \{N(x): x \in \operatorname{cl}(M \setminus \mathcal X_T)\}$. Since $\operatorname{cl}(M \setminus \mathcal X_T)$ is compact, $N^*$ is finite.
	
	Let $x \in M$. Therefore there exists a control sequence $\{u_{k+i|k}\}_{k \in \mathbb Z_+}$ guiding $x_k = x$ 
	to $x_e$ and $x_{k+N^*|k} \in \mathcal X_T$.
	The sequence $\{u_{k+i|k}\}_{k \in \mathbb Z_{N^*}}$ is feasible for (\ref{equ:mpcopt}) because it satisfies both state and control constraints. Therefore, $\mathcal D_{N^*} \supset \mathcal D_\infty = M$.
	
	By Theorem 1, $M$ is the domain of attraction of the closed-loop system defined by (\ref{equ:dyn}) and (\ref{equ:mpclaw}), and therefore the closed-loop system is globally asymptotically stable at $x_e$.
	\qed
	
	


	\section{Local control law}\label{sec:locallaw}
	
	The development of the control law (\ref{equ:mpclaw}) depends on the design of target set $\mathcal X_T$ and the local control law $\kappa: M \rightarrow U$ with the properties enumerated in (\ref{equ:locallaw}). 
	In this section, we describe a procedure to construct $\mathcal X_T$ and $\kappa$ 
	that satisfy the properties in \eqref{equ:locallaw} for dynamics \eqref{equ:dyn} evolving on manifolds.
	For conventional  MPC, the terminal set and terminal cost are designed from a linearization of the dynamics around the equilibrium \cite{rawlings_mpc_book}. In our approach, 
	we utilize diffeomorphisms to obtain local coordinates, on which we 
	construct a stabilizing controller; we then transform the local control law to the coordinates of the original system.
	
	
	
	To begin, because $M$ is an $n$-dimensional manifold and $U$ is a compact subset of an $m$-dimensional manifold that contains $u_e$ in its interior, there exist local diffeomorphisms at $x_e$ and $u_e$ \cite{milnor_book,guillemin_book},
	\begin{equation}
	\phi: V \rightarrow M, \psi: W \rightarrow U,
	\end{equation}
	where $V \subset \mathbb R^n$ and $W \subset \mathbb R^m$ are open neighborhoods of $x_e$ and $u_e$, respectively. 
	Furthermore, because the equilibrium is in the interior of both $V$ and $W$, $f$ is continuously differentiable, and $f(x_e,u_e) = x_e$, there exists
	an open neighborhood $V' \times W' \subset V \times W$ of $(0,0)$ such that $f(\phi(V') \times \psi(W')) \subset \phi(V)$. Therefore,
	\begin{equation}
	f' = \phi^{-1} \circ f \circ (\phi \times \psi),
	\end{equation}
	and the derivative of $f'$ at $(0,0)$ is,
	\begin{equation}
	df'_{(0,0)} = d\phi_0^{-1} \circ df_{(x_e,u_e)} \circ (d\phi_0 \times d\psi_0).
	\end{equation}
	
	Let $\xi_k \in \mathbb R^n$ and $\eta_k \in \mathbb R^m$ and 
	let $A = df'_{(0,0)} \circ (I_n,0_{m \times n})$ and $B = df'_{(0,0)} \circ (0_{n \times m},I_m)$. Define a linear update equation,
	\begin{equation}\label{equ:lindyn}
	\xi_{k+1} = A\xi_k + B\eta_k.
	\end{equation}
	
	Let,
	\begin{equation}\label{equ:Lprime}
	L' = \frac{1}{\lambda}L \circ (\phi \times \psi),
	\end{equation}
	where $0 < \lambda < 1$ is a scalar parameter. 
	Consider the discrete-time algebraic Riccati equation,
	\begin{align}
	&0 = A^\textrm{T}PA - P + Q \label{equ:dare} \\ &- (A^\textrm{T}PB+N)(BPB^\textrm{T}+R)^{-1}(A^\textrm{T}PB+N)^\textrm{T}, \nonumber
	\end{align}
	where,
	\begin{equation}\label{equ:Hess}
	\begin{bmatrix}Q & N \\ N^\textrm T & R\end{bmatrix} = \operatorname{Hess}L'(0,0).
	\end{equation}
	is the Hessian of $L'$ at $(0,0)$. A positive-definite solution $P$ to the algebraic Riccati equation (\ref{equ:dare}) exists if and only if the pair $(A,B)$ is stabilizable and the associated symplectic pencil does not have eigenvalues on the unit circle \cite{arnoldlaub84}. A sufficient condition implying the latter is that $\operatorname{Hess}L'(0,0)$ is positive-definite.
	
	Suppose a positive-definite solution $P$ exists for (\ref{equ:dare}). Let $\kappa': \mathbb R^n \rightarrow \mathbb R^m$ be a stabilizing feedback control law for (\ref{equ:lindyn}) where,
	\begin{equation}\label{equ:linlaw}
	\kappa'(\xi_k) = -(BPB^\textrm{T}+R)^{-1}(A^\textrm{T}PB+N)^\textrm{T}\xi_k.
	\end{equation}
	We introduce a set $\mathcal P_c \subset \mathbb R^n$ and a function $F': \mathbb R^n \rightarrow \mathbb R$ where,
	\begin{align}
	\mathcal P_c &= \{\xi \in \mathbb R^n: F'(\xi) \leq c\}, \label{equ:Pcdef} \\
	F'(\xi) &= \xi^\textrm TP\xi. \label{equ:Fprime}
	\end{align}
	Because $P$ is the solution to the algebraic Riccati equation (\ref{equ:dare}), the design of the control law (\ref{equ:linlaw}) implies that the set $\mathcal P_c$ is compact and invariant with respect to the closed loop dynamics (\ref{equ:lindyn}), (\ref{equ:linlaw}). 
	
	Our goal is to use $\mathcal P_c$ and $\kappa'$ in order to design the target set $\mathcal X_T$ and local control law $\kappa$. The following result is a proof of the existence of $\mathcal X_T$ and $\kappa$.
	\begin{prop}
		Let $F = F' \circ \phi^{-1}$. Suppose there exists a solution $P > 0$ to the algebraic Riccati equation (\ref{equ:dare}). Then there exists $c > 0$ such that the set $\mathcal X_T \subset \phi(\mathcal P_c)$ and control law $\kappa = \kappa' \circ \phi^{-1}$ satisfy the assumptions of (\ref{equ:locallaw}).
	\end{prop}
	\textit{Proof.}
	Firstly, we show that $F$ satisfies the assumptions made in the design of the MPC law. Because $M$ is a smooth manifold, $\phi$ is smooth and this, along with the smoothness of $F'$, implies that $F$ is smooth. Furthermore, $F(x_e) = F'(\phi^{-1}(x_e)) = F'(0) = 0$ which satisfies (\ref{equ:costprops1}). Finally, due to the equivalence of norms on $\mathbb R^n$ and (\ref{equ:costprops2}), there exists $r > 0$ such that $F(x) = F'(\phi^{-1}(x)) = \phi^{-1}(x)^{\textrm{T}}P\phi^{-1}(x) \geq r\phi^{-1}(x)^{\textrm{T}}Q\phi^{-1}(x) = rL'(\phi^{-1}(x),0) = rL(x,u_e) \geq r\cdot\gamma(d(x,x_e))$. Since $r\cdot\gamma$ is a class $\mathcal K$ function, (\ref{equ:costprops3}) is satisfied.
	
	If $c > 0$, then by definition (\ref{equ:Pcdef}), the set $\mathcal P_c$ has a non-empty interior, which contains the origin. Because $V'$ contains the origin, there exists $c' > 0$ such that $\mathcal P_{c'} \subset V'$ and $\kappa'(\mathcal P_{c'}) \subset U'$ and therefore $\phi(\mathcal P_{c'}) = \{x: F(x) \leq c'\}$ because $F'(\phi^{-1}(x)) = F(x)$.
	
	Furthermore, because $(x_e,u_e)$ is in the interior of $\mathcal X \times \mathcal U$ and $\kappa'$ is $C^1$, there exists $c'' > 0$ such that $(\phi(\xi),\kappa(\phi(\xi))) \in \mathcal X \times \mathcal U$ for all $\xi \in \mathcal P_{c''}$. 
	It remains to show that there exists a $c > 0$ such that $\mathcal X_T = \phi(\mathcal P_c)$ satisfies conditions (\ref{equ:locallaw2}) and (\ref{equ:locallaw3}). 
	
	
	
	Note that,
	\begin{align*}
	f'(\xi_k,\eta_k) &= A\xi_k+B\eta_k+o(\|(\xi_k,\eta_k)\|), \\
	L'(\xi_k,\eta_k) &= \begin{bmatrix}\xi_k^{\textrm T} & \eta_k^{\textrm T}\end{bmatrix}\begin{bmatrix}Q & N \\ N^\textrm T & R\end{bmatrix}\begin{bmatrix}\xi_k \\ \eta_k\end{bmatrix} \\ &\hspace{1in} +o(\|(\xi_k,\eta_k)\|^2), \\
	F'(\xi_k) &= \xi_k^\textrm T P\xi_k.
	\end{align*}
	Therefore, $F'(\xi_{k+1})-F'(\xi_k)+ L'(\xi_k,\kappa'(\xi_k)) = F' \circ f'(\xi_k,\kappa'(\xi_k))-F'(\xi_k)+ L'(\xi_k,\kappa'(\xi_k)) = o(\|\xi_k\|^2)$.
	This implies that for $x_k \in \mathcal P_{c'}$,
	$F(x_{k+1})-F(x_k)+ L(x_k,\kappa(x_k)) 
	= F \circ f(x_k,\kappa(x_k))-F(x_k)+ L(x_k,\kappa(x_k))
	= F' \circ f' (\phi^{-1}(x_k),\kappa(\phi^{-1}(x_k)))-F'(\phi^{-1}(x_k))
	+ \lambda L'(\phi^{-1}(x_k),\kappa(\phi^{-1}(x_k))) + o(\|\phi^{-1}(x_k)\|^2)$, and therefore,
	\begin{align*}
	&F(x_{k+1})-F(x_k)+ L(x_k,\kappa(x_k)) \\ &~= -(1-\lambda) L'(\phi^{-1}(x_k),\kappa(\phi^{-1}(x_k))) + o(\|x_k\|^2).
	\end{align*}
	Because $L'$ is positive definite, for any choice of $0 < \lambda < 1$, there exists $c''' \geq 0$ such that  (\ref{equ:locallaw2}) and (\ref{equ:locallaw3}) are satisfied. 
	Let $c = \min(c'',c''')$ to complete the proof.
	\qed

	\section{Globally stabilizing control law}\label{sec:globallaw}
	
	The results of Section \ref{sec:mpc} show that the MPC law (\ref{equ:mpclaw}) is globally asymptotically stabilizing under certain assumptions. We will now show that this implies that the control law generated by MPC is necessarily discontinuous for certain classes of manifolds. Specifically, we will show that there does not exist a globally stabilizing, continuously differentiable control law for compact manifolds $M$ with Euler characteristic $\chi(M)$ not equal to 1. The Euler characteristic is a topological invariant that restricts the possible combinations of sinks and sources that can exist for a continuous vector field defined on the manifold \cite{guillemin_book}. Generally, only a manifold with Euler characteristic of 1 admits a continuous vector field with one sole sink and no other equilibrium. As an example, the normed vector space $\mathbb R^n$ has an Euler characteristic of 1, the sphere $S^2$ has an Euler characteristic of 2, and all compact Lie groups have an Euler characteristic of 0. Therefore, the sphere and compact Lie groups do not admit continuous vector fields with one sink and no other equilibrium.
	
	We begin by noting that 
	there exists no globally stabilizing, continuously differentiable, \textit{continuous-time} control law on compact manifolds $M$  \cite{bhatbernstein00}. 
	The proof of the continuous-time result appeals to a consequence of the Poincar\'{e}-Hopf theorem, which states that the sum of fixed point indexes on $M$ must be equal to the Euler characteristic. Since the index of a sink is equal to 1, if $\chi(M) \neq 1$, the set of fixed points on $M$ cannot consist of solely a sink. In discrete-time, the analogous result to the Poincar\' e-Hopf theorem is the Lefschetz-Hopf theorem, of which we use a strong version below.
	
	\begin{thm}[Lefschetz-Hopf \cite{granas_book}]
		Let $M$ be a compact smooth manifold and let $f_0: M \rightarrow M$ be a continuous map that is homotopic to the identity map $\textrm{id}: M \rightarrow M$. Suppose the set of points $x \in M$ satisfying $x = f_0(x)$ is finite. Then, 
		\begin{equation}
		\sum_{x_0 \in \{x: x = f_0(x)\}} i_{f_0}(x_0) = \chi(M),
		\end{equation}
		where $i_{f_0}(x_0)$ is the fixed point index of $f_0$ at $x_0$.
	\end{thm}
	
The fixed point index is a property characterizing an equilibrium. 
Note that when the fixed point is an asymptotically stable equilibrium, \textit{i.e.}, a sink, then its fixed point index is equal to 1. This fact is used in the proof of the discrete-time analogue to Brockett's necessary condition \cite{linbyrnes94} for systems where $M = \mathbb R^n$ but,
to the authors' knowledge, an equivalent result for $M$ with $\chi(M) \neq 1$ does not appear elsewhere in the literature. Hence we present below a theorem proving that there exists no globally-stabilizing, continuous, \textit{discrete-time} control law on manifolds $M$ with $\chi(M) \neq 1$ whenever the closed-loop dynamics are homotopic to the identity.
	
	\begin{thm}\label{thm:c1}
		Let $M$ be a compact smooth manifold. Consider the dynamics (\ref{equ:dyn}) and assume the unforced dynamics map $f_0: x \mapsto f(x,u_e)$ is homotopic to identity. Further assume that there exists a 
		control law $u_k = \kappa(x_k)$ with globally asymptotically stable equilibrium $x = x_e$, and define the resulting closed-loop dynamics by $x_{k+1} = f_1(x_k) = f(x_k,\kappa(x_k))$. If $\chi(M) \neq 1$, then $\kappa: M \rightarrow U$ is not continuous.
	\end{thm}
	
	\begin{rem}
		The homotopy assumption is crucial since it ensures that the unforced dynamics correspond to a continuous deformation of the system from an equilibrium configuration. A case where this assumption is violated is when the unforced map satisfies $f_0(x_k) = x_e$ for all $x_k \in M$. In this case, $\kappa(x_k) = 0$ is an asymptotically stabilizing continuous control law. Nevertheless, for mechanical systems, the assumption is generally valid because, for such systems, 
		the dynamics (\ref{equ:dyn}) correspond to a discretization of a continuous-time system, \textit{i.e.}, if $x(0) = x_k$ and $u_k = 0$, then $x_{k+1} = \phi(x(0),h)$, where $h$ is the discretization step and $\phi$ is the continuous flow map. 
	\end{rem}
	
	\begin{cor}
		As a result of Proposition \ref{prop:gas} and Theorem \ref{thm:c1}, if there exists a control sequence that is 
		capable of steering any initial condition to $x_e$, then there exists a finite $N$ such that the MPC law is globally asymptotically stabilizing and therefore the MPC law is discontinuous.
	\end{cor}
	
	\textit{Proof of Theorem \ref{thm:c1}.}
	Assume $\kappa$ is continuous. This implies that $f_1$ is continuous. 
	By assumption, $\lim_{k \rightarrow \infty} f_1(x_k) = x_e$ for all $x_0 \in M$. This fact, along with the continuity of $f_1$ implies that $M$ is connected.
	All connected manifolds are also path-connected, so $M$ is path-connected. Define a continuous path $p_{x_k}: [0,1] \rightarrow M$ from $x_e$ to $x_k$. Since $f_0$ is homotopic to the identity, there exists a homotopy $H_0$ from $\operatorname{id}_M$ to $f_0$. Define the map $H_1: M \times [0,1] \rightarrow M$ such that $H_1(x_k,t) = H_0(x_k,2t)$ for $t \in [0,\frac{1}{2})$ and $H_1(x_k,t) = f(x_k,\kappa(p_{x_k}(2t-1)))$ for $t \in [\frac{1}{2},1]$. Then, $H_1(x_k,0) = x_k$ and $H_1(x_k,1) = f_1(x_k)$. Because $\kappa(x_e) = u_e$, $H_1$ is continuous and is therefore a homotopy between $\operatorname{id}_M$ and $f_1$, \textit{i.e.}, $f_1$ is homotopic to the identity.
	
	The sole equilibrium is $x_e$ and $i_{f_1}(x_e) = 1$. Therefore, $\sum_{x_0 \in \{x: x = {f_1}(x)\}} i_{f_1}(x_0) = i_{f_1}(x_e) = 1$. However, $\chi(M) \neq 1$, and, as a consequence of the Lefschetz-Hopf theorem, this leads to a contradiction. 
	Therefore, $\kappa$ cannot be continuous. 
	\qed

	\section{Application to $\textrm{SO}(3)$ and spacecraft attitude control}\label{sec:so3}
	
	In this section, we apply the previously developed techniques to the constrained control of spacecraft attitude.\footnote{This example originally appeared in \cite{kalabicgupta_acc14} and in this paper, we report additional details, derivations, and include further discussion relevant to this example} 
	The orientation of a spacecraft can be uniquely represented by an element of the Lie group $\textrm{SO}(3)$, which is the group of all orthogonal matrices whose determinant is equal to 1.
	Physically, the first, second, and third columns of an element $g \in \textrm{SO}(3)$ represent the direction of the $x$, $y$, and $z$ axes, respectively, when viewed from a fixed frame in $\mathbb R^3$. This 
	allows us to derive the discrete-time dynamics using Lie group variational integrator techniques, resulting in a dynamic equation that ensures that state updates are elements of $\textrm{SO}(3)$ within a guaranteed numerical tolerance.
	The LGVI spacecraft dynamics are given by \eqref{equ:so3_dyn_intro}.
	In \eqref{equ:so3_dyn_intro}, $g_k \in \textrm{SO}(3)$ represents the spacecraft orientation, $f_k \in \textrm{SO}(3)$ is a one time-step change in $g_k \in \textrm{SO}(3)$, 
	and $u_k \in \mathfrak{so}(3)$ is related to the applied torque $\tau_k$ by the equation $u_k = \tau_k^{\times}$, where $\mathfrak{so}(3)$ is the set of skew-symmetric 3-by-3 matrices\footnote{$\cdot^\times$ is the map from elements of $\mathbb R^3$ to elements of $\mathfrak{so}(3)$ which preserves the cross product under multiplication, \textit{i.e.}, if $C = a^\times$, then $a \times b = Cb$ for any $b \in \mathbb R^3$, where $\times$ is the cross-product.}. 
	Note that the dynamics (1) of \cite{kalabicgupta_acc14} are equivalent to the dynamics (\ref{equ:so3_dyn_intro}) in this paper. In \cite{kalabicgupta_acc14}, the dynamics are presented in Hamiltonian form, whereas here we use the Lagrangian form (see \cite{lee_diss}, which presents dynamics in both forms, for details).

	To solve the implicit equation (\ref{equ:so3_fdyn_intro}), we use the procedure from \cite{cardosleite03}. At time $k$, the quantities $f_k$ and $u_k$ are known, so we let $M_k = Jf_k-f_k^{\textrm T}J+h^2u_k$ so that $f_{k+1}J-Jf_{k+1}^{\textrm T} = M_k$. Then,
	\begin{equation}
	f_{k+1} = \left(M_k/2+S_k\right)J^{-1},
	\end{equation}
	where $S$ is the solution to the algebraic Riccati equation,
	\begin{equation}\label{equ:so3_are}
	-(M_k/2)^\textrm{T}S_k-S_k(M_k/2)-S_k^2+(J^2+M_k^2/4) = 0.
	\end{equation}
	Note that (\ref{equ:so3_are}) is solvable if and only if the term $J^2+M_k^2/4$ is positive semi-definite, \textit{i.e.},
	\begin{equation}
	J^2+M_k^2/4 \succeq 0,
	\end{equation}
	which is 
	a convex constraint that can be enforced by the MPC law \cite{kalabicgupta_acc14}.
	
	
	
	\subsection{MPC law}
	
	In this section, we develop an MPC-based 
	spacecraft attitude controller. We begin by choosing an appropriate cost function of the form (\ref{equ:mpccost}) with,
	\begin{multline}\label{equ:so3_cost}
	L(g_k,f_k,u_k) = \operatorname{tr}(Q_g(I_3-g_k)) \\ +\frac{1}{h^2}\operatorname{tr}(Q_f(I_3-f_k)) +\frac{1}{2}\operatorname{tr}(u_k^{\textrm T}Ru_k),
	\end{multline}
	and with positive-definite symmetric matrices $Q_g$, $Q_f$, and $R$.
	Next we construct a locally stabilizing control law.
	
	The linearized dynamics (\ref{equ:lindyn}) corresponding to (\ref{equ:so3_dyn_intro}) evolve on $\mathbb R^3 \times \mathbb R^3$ and, according to \cite{leeleok08}, are given by, 
	\begin{equation}\label{equ:linearization}
	\begin{bmatrix}\zeta_{k+1} \\ \omega_{k+1}\end{bmatrix} = A\begin{bmatrix}\zeta_{k} \\ \omega_{k}\end{bmatrix} + B\tau_k,
	\end{equation}
	where $\zeta_k, \omega_k \in \mathbb R^3$ satisfy $\zeta_k^\times = \operatorname{Log} g_k$ and $(h\omega_k)^\times = \operatorname{Log} f_k$, where $\operatorname{Log}$ is the standard branch of the matrix logarithm function. 
	The matrices $A$ and $B$ are,
	\begin{equation}
	A = \begin{bmatrix}I_3 & h \\ 0 & I_3\end{bmatrix},~B = \begin{bmatrix}0 \\ h I_3\end{bmatrix}.
	\end{equation}
	In the language of Section \ref{sec:locallaw}, the composition of the map $\cdot^\times$ and the matrix exponential $\exp$ from $\zeta_k$ to $g_k$ and from $h\omega_k$ to $f_k$ give the diffeomorphism $\phi$. The diffeomorphism $\psi$ is the map $\cdot^\times$. We now define $L'$ according to (\ref{equ:Lprime}). Specifically, we choose $0 < \lambda < 1$ so that,
	\begin{multline*}
	\lambda L'(\zeta_k,\omega_k,\tau_k) = 
	\operatorname{tr}(Q_g(I_3-\exp(\zeta_k^\times))) \\ +\frac{1}{h^2}\operatorname{tr}(Q_f(I_3-\exp(h\omega_k)^\times))+\frac{1}{2}\operatorname{tr}(\tau_k^{\times\textrm T}R\tau_k^{\times}), 
	\end{multline*}
	which implies that,
	\begin{multline}
	\lambda L'(\zeta_k,\omega_k,\tau_k) = 
	\frac{1}{2}\zeta_k^{\textrm{T}}\tilde Q_g\zeta_k + \frac{1}{2}\omega_k^{\textrm{T}}\tilde Q_f\omega_k + \frac{1}{2}\tau_k^{\textrm{T}}\tilde R\tau_k
	\\ + o(\|(\zeta_k,\omega_k)\|^2), \label{equ:so3hess}
	\end{multline}
	where $\tilde {Q}_{g} = \operatorname{tr}({ Q}_{g})I_3 - { Q}_{g}$, 
	$\tilde { Q}_{f} = \operatorname{tr}({Q}_{f})I_3 - {Q}_{f}$, 
	and $\tilde R = \operatorname{tr}(R)I_3 - R$. 
	
	The equation (\ref{equ:so3hess}) can be used to define the Hessian matrix of (\ref{equ:Hess}). We use (\ref{equ:so3hess}) in order to construct the terminal cost $F$ and target set $\mathcal X_T \subset \mathcal P_c$ according to the steps outlined in (\ref{equ:dare})-(\ref{equ:Fprime}). 
	
	Thus, the MPC law (\ref{equ:mpclaw}) has been defined by the finite horizon optimal control problem \eqref{equ:mpcopt}, which has been formulated by using the constraint sets $\mathcal X$ and $\mathcal U$ together with the cost functions $F$ and $L$ and the terminal set $\mathcal X_T$ designed above. 
	
	We now consider rest-to-rest controllability for the dynamics (\ref{equ:so3_dyn_intro}) subject to constraints. 
	We show that, given a long enough prediction horizon, the MPC control law is able to stabilize to the origin any spacecraft attitude $g_k$ that is near to rest.
	In the case that we consider, there is no constraint on the spacecraft rotation $g_k$; however there may be limitations on the rate of change $f_k$ or allowed torque $\tau_k$.
	
	\begin{cor}
		Let the set $\mathcal D_N$ be defined according to the definition in (\ref{equ:Dn_defn}) and the dynamics in (\ref{equ:so3_dyn_intro}).
		Assume $\textrm{SO}(3) \times \{I_3\} \subset \mathcal X$, where $\mathcal X$ satisfies the properties given in Section \ref{sec:mpc}. Then there exists a finite number $N$ such that $\textrm{SO}(3) \times V' \subset \mathcal D_N$, where $V' \subset \textrm{SO}(3)$ is an open set containing $I_3$.
	\end{cor}
	
	The proof of the corollary is available in \cite{kalabicgupta_acc14}.

	\subsection{Simulation results}
	
	We present a numerical simulation of the MPC law on $\textrm{SO}(3)$. 
	In the simulation, in order to provide an example of discontinuity in the MPC law, we simulate a rest-to-rest rotation from a rotation of 180 degrees to $I_3$.
	
	For the cost function (\ref{equ:so3_cost}), we choose the following matrices: $Q_g = I_3$, $Q_f = J$, $R = 2I_3$, and $\lambda = 0.1$. In both cases, the terminal constraint set is chosen as $\mathcal X_T = \mathcal P_c$, where $\mathcal P_c$ is defined as in (\ref{equ:Pcdef}) and is chosen to be as large as possible. We 
	set $N = 10$ 
	and, $\mathcal U$ is large enough that the constraints are not active during our simulation.
	
	
	%
	%
	
\begin{figure}[p]
		\centering
		\includegraphics[width = 0.45\textwidth]{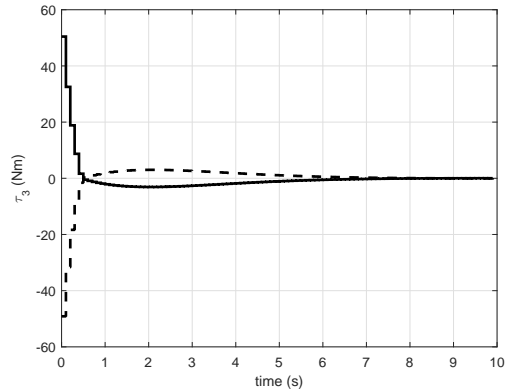}
		\caption{$\tau_{k,3}$ for initial rotation of 180 degrees (solid) and $-0.99 \cdot 180$ (dotted)}
		\label{fig:torque}
	\end{figure}

\begin{figure}[p]
		\centering
		\includegraphics[width = 0.45\textwidth]{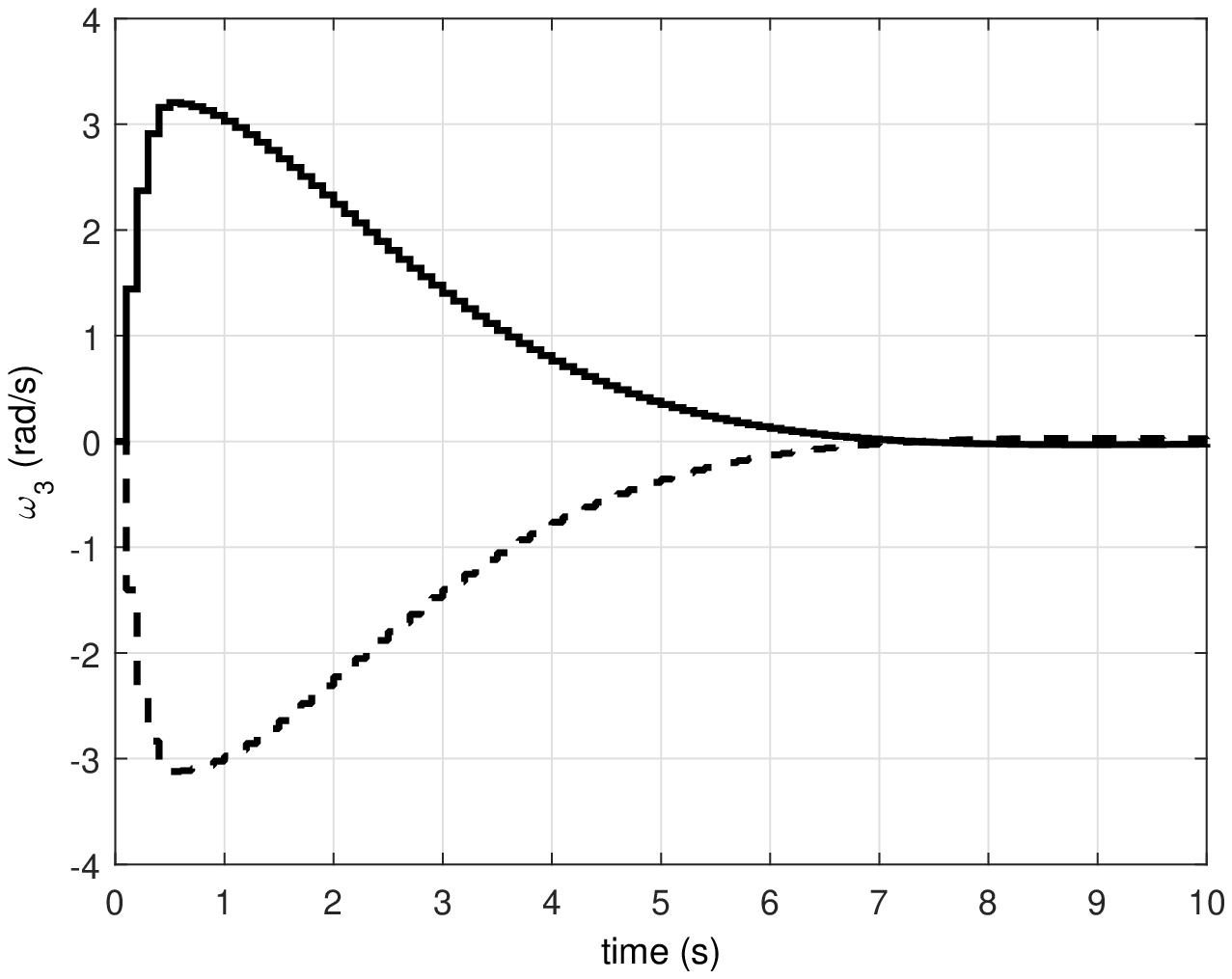}
		\caption{$\omega_{k,3}$ for initial rotation of 180 degrees (solid) and $-0.99 \cdot 180$ (dotted)}
		\label{fig:omega}
	\end{figure}

\begin{figure}[p]
		\centering
		\includegraphics[width = 0.45\textwidth]{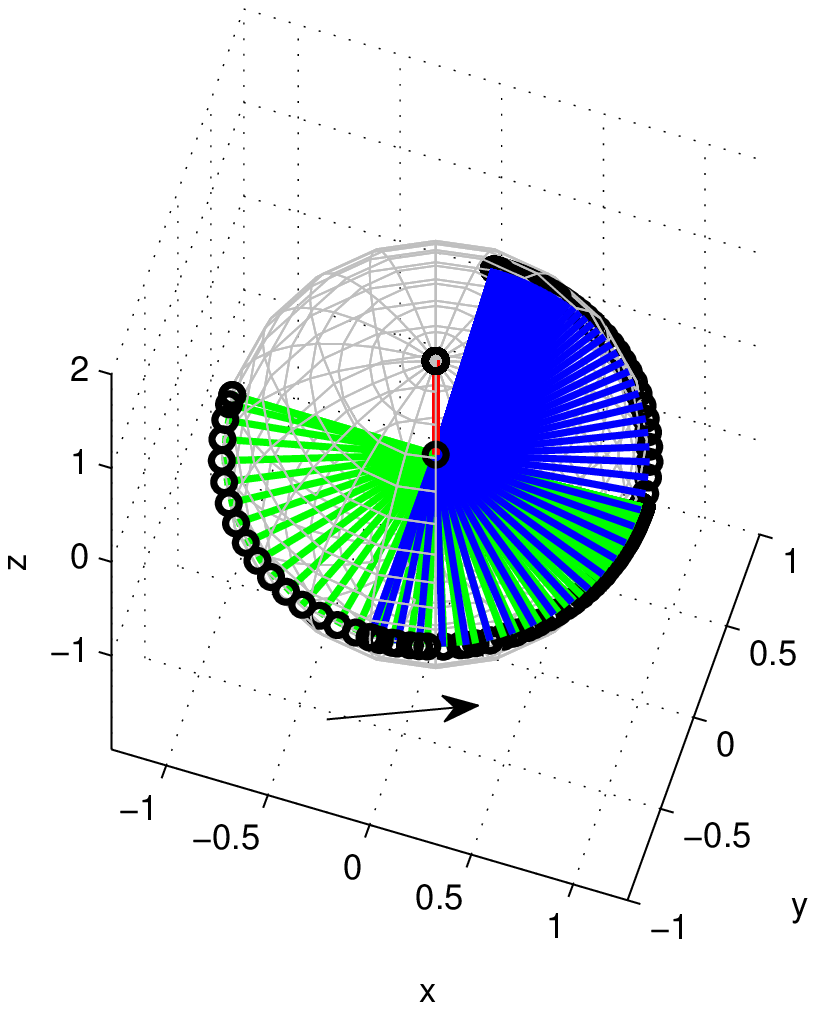}
		\caption{Orientation maneuver corresponding to the initial rotation of 180 degrees plotted at 2s increments}
		\label{fig:ram_180}
	\end{figure}

\begin{figure}[t]
		\centering
		\includegraphics[width = 0.45\textwidth]{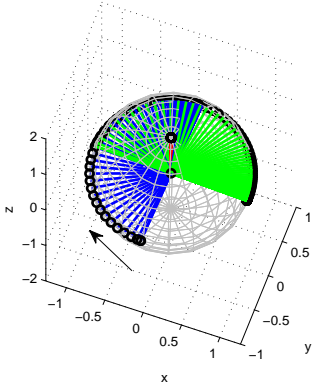}
		\caption{Orientation maneuver corresponding to the initial rotation of $-0.99 \cdot 180$ degrees plotted at 2s increments}
		\label{fig:ram_179}
	\end{figure}

	Theoretically, a control law that is globally stabilizing should exhibit a discontinuity, so our goal for the simulation 
	is to confirm that this is true of our controller. 
	The results for two simulations are presented in Figs.~\ref{fig:torque}-\ref{fig:ram_179}. In the first simulation, the initial rotation is 180 degrees about the $z$-axis, and in the second case, it is close but not equal to $-180$ degrees about the $z$-axis. The torque and angular velocity on the $z$-axis are shown in Figs.~\ref{fig:torque}-\ref{fig:omega}, respectively, where we can observe that, although the initial conditions are very close to each other, the trajectories are almost opposite in sign. This shows that there is a discontinuity in the control law. In fact, the discontinuity occurs exactly at the location of the 180 degree rotation or at any resting initial condition $g_0$ for which $\operatorname{tr}(g_0) = -1$, because these points lie on the branch cut of the $\operatorname{Log}$ function. 
	
	In Figs.~\ref{fig:ram_180}-\ref{fig:ram_179}, the direction of rotation is marked with an arrow. We see that in one case, the control rotates the spacecraft counter-clockwise whereas, in the other case, the rotation is clockwise.
	
	\section{Conclusion}\label{sec:conc}
	
	This paper presented a general MPC theory for dynamics that evolve on smooth manifolds. An MPC scheme was developed for manifolds that generalizes the results of conventional MPC in the case of $\mathbb R^n$. In the case of MPC on manifolds, a general construction for a locally valid control law on manifolds was derived, which was used in the design of the MPC terminal cost and target set.
	
	Results were presented that showed the globally stabilizing properties of the MPC scheme as well as the nonexistence of continuous discrete-time control laws on manifolds whose Euler characteristic is not equal to 1. 
	Thus, MPC achieves stability by producing a discontinuous control law, even in cases where a continuous globally stabilizing law does not exist. 
	The authors believe that the ability to generate, possibly discontinuous, stabilizing feedback control laws for systems whose dynamics evolve on manifold is appealing for practical applications.  
	
	
	Finally, 
	an application of the results to a constrained spacecraft attitude control problem was presented in the case of the matrix Lie group $\textrm{SO}(3)$. The simulation results that were presented showed the expected discontinuity in the MPC law.
	
	\appendix

	\section{Proof of Theorem \ref{thm:mpc}}
	
	In order to prove Theorem \ref{thm:mpc}, we rely on discrete-time Lyapunov stability theory for smooth manifolds. Therefore, we begin by presenting the Lyapunov stability analysis. Let $M$ be an $n$-dimensional smooth manifold with metric $d$. Consider the following discrete-time dynamical system,
	\begin{equation}\label{A1}
	x_{k+1} = f(x_{k}),
	\end{equation}
	where $x_k \in M$ and $f: M\to M$ is a continuous function. 
	
	\begin{definition}
		A point $x_{e}\in M$ is called an equilibrium point of (\ref{A1}) if $f(x_{e}) = x_{e}$.
	\end{definition}
	
	Note that if $f(x_{e}) = x_{e}$, then $f^{k}(x_{e}) = x_{e}$, for all $k \in \mathbb Z_+$, where $f^{k+1} = f \circ f^{k}$ for all $k \in \mathbb Z_+$ and $f^0 = \operatorname{id}_M$.
	
	\begin{definition}
		An equilibrium point $x_{e}\in M$ is said to be Lyapunov stable if for any open neighborhood ${U} \subset M$ of $x_{e}$, there exists an open neighborhood ${W} \subset M$ of $x_{e}$ such that for all $x_{0}\in{W}$, $f^{k}(x_{0})\in{U}$, for all $k \in \mathbb Z_+$. It is said to be locally asymptotically stable if it is Lyapunov stable and for all $x_0 \in W$, $\lim_{k \to \infty} f^{k}(x_{0}) = x_{e}$ or, equivalently, $\lim_{k \to \infty} d(f^{k}(x_{0}),x_{e}) = 0$.
	\end{definition}
	
	
	\begin{definition} 
		Let $x_{e}$ be an equilibrium point of (\ref{A1}). A function $\mathcal V: M\to\mathbb{R}$ is a Lyapunov function if there exist two class-$\mathcal K$ functions $\alpha_1,\alpha_2$, that satisfy $\alpha_1(d_e(x)) \leq \mathcal V(x) \leq \alpha_2(d_e(x))$ for all $x$ in some neighborhood of $x_e$,
		and $\Delta \mathcal V(x) \colonequals \mathcal V(f(x))-\mathcal V(x)$ is 
		negative semi-definite in a neighborhood of $x_{e}$. If $\Delta \mathcal V(x)$ is 
		negative definite in a neighborhood of $x_e$, 
		then $\mathcal V$ is a strict Lyapunov function.
	\end{definition}
	
	More formally, the above definition implies that $\mathcal V$ is a Lyapunov function if there exists a neighborhood ${V}$ of $x_{e}$ such that $\mathcal V(x)>0$, $\Delta \mathcal V(x)\leq 0$ for all $x\in{V}\setminus\{x_{e}\}$, and $\mathcal V(x_{e}) = 0$. If, instead of $\Delta \mathcal V(x) \leq 0$, the condition $\Delta \mathcal V(x)<0$ holds for all $x\in{V}\setminus\{x_{e}\}$, then the Lyapunov function is a strict Lyapunov function.
		
	\begin{thm}
		Let $x_{e}$ be an equilibrium point of (\ref{A1}).
		\begin{enumerate}[(i)]
			\item{If there exists a Lyapunov function $\mathcal V$, then $x_{e}$ is Lyapunov stable.}
			\item{If there exists a strict Lyapunov function $\mathcal V$, then $x_{e}$ is locally asymptotically stable.}
		\end{enumerate}
	\end{thm}
	
	\textit{Proof.} The proof follows arguments similar to the one given for the continuous-time case in \cite{bullo_book}. 
	\begin{enumerate}[(i)]
		\item 
		Let $U \subset M$ be a neighborhood of $x_e$. 
		As a consequence of Lemma 2.23 in \cite{lee_book}, $M$ admits a locally finite cover by precompact sets, implying that there exists a neighborhood $U' \subset M$ of $x_e$ that is a subset of the union of finitely many precompact sets, and therefore is itself precompact.
		Let $U'' \subset U \cap U'$ be an open neighborhood of $x_e$ on which $\alpha_1(d_e(x)) \leq \mathcal V(x) \leq \alpha_2(d_e(x))$ for all $x \in U''$.
		Let $V_\varepsilon := \{x \in M: d_e(x)\leq \varepsilon\}$ and choose $\varepsilon > 0$ so that $V_\varepsilon$ is a subset of $U''$.
		By construction, 
		$V_\varepsilon$
		is a closed subset of a precompact set and is therefore compact. 
		Let $V_\varepsilon(x_e) \subset V_\varepsilon$ be the connected component of $V_\varepsilon$ containing $x_e$.
		The interior of $V_\varepsilon(x_e)$ is non-empty and contains $x_e$. 
		Let $x_0\in\operatorname{int}V_\delta(x_{e})$ where $\delta := \alpha_2^{-1}(\alpha_1(\varepsilon)) \leq \varepsilon$. Since $\mathcal V(x_0)\leq\alpha_2(d_e(x_0)) \leq \alpha_2(\delta)$ and $\mathcal V(f^k(x_0))-\mathcal V(x_0) = \sum_{j=1}^{k}(\mathcal V(f^j(x_0))-\mathcal V(f^{j-1}(x_0)))\leq 0$, it follows that $\mathcal V(f^k(x_0)) \leq \mathcal V(x_0) \leq \alpha_2(\delta)$, which implies that 
		$d_e(f^k(x_0)) \leq \alpha_1^{-1}(\mathcal V(f^k(x_0))) \leq \alpha_1^{-1}(\alpha_1(\varepsilon)) = \varepsilon$ for all $k \in \mathbb Z_+$.
		Because $V_\varepsilon$ is compact and $\mathcal V$ is positive definite in a neighborhood of $x_e$, there exists a constant $\varepsilon' \leq \varepsilon$ 
		such that $V_{\varepsilon'}$ consists of only one connected component, so that $V_{\varepsilon'} = V_{\varepsilon'}(x_e)$. As a consequence, for any neighborhood $U$, there exists $W = \operatorname{int}V_{\varepsilon'}$ such that for any $x_0 \in W$, $f^{k}(x_0) \in U$ for all $k \in \mathbb Z_+$. \qed
		
		\item
		Let $U \subset M$ be a neighborhood of $x_e$. Choose $\varepsilon' > 0$ so that $\mathcal V$ is positive definite, 
		$\Delta \mathcal V$ is strictly negative definite on $V_{\varepsilon'}$, and $V_{\varepsilon'}$ is a compact set with only one connected component. Let $W = \operatorname{int}V_{\varepsilon'}$ as in the proof of (i), implying that $x_e$ is Lyapunov stable. For all $x_0 \in W$, the sequence $\{\alpha_1(d_e(f^{k}(x_0)))\}_{k \in \mathbb Z_+}$ is a non-increasing sequence bounded below by $0$, and hence $\lim_{k\to\infty} \alpha_1(d_e(f^{k}(x_0))) = \beta \geq 0$. We can show that $\beta = 0$ by contradiction. Assume $\beta > 0$, implying that the set $V_{\varepsilon'} \setminus \operatorname{int}V_{\alpha_1^{-1}(\beta)}$ does not contain $x_e$. Because this set is compact and does not contain $x_e$, the strictly negative definite, continuous function $\Delta \mathcal V$ attains a negative maximum on this set. Denote the maximum by $-\sigma < 0$. For any $x_0 \in W$, $\mathcal V(f^{k}(x_0)) = \mathcal V(x_{0})+\sum_{j=1}^{k}(\mathcal V(f^{j}(x_{0}))-\mathcal V(f^{j-1}(x_{0})))\leq \mathcal V(x_{0})-k\sigma$. For any $k > \mathcal V(x_{0})/\sigma$, it follows that $\alpha_1(d_e(f^{k}(x_{0}))) \leq \mathcal V(f^{k}(x_{0}))<0$, which contradicts the fact that $\alpha_1\circ d_e$ is positive definite. Therefore $\beta = 0$ and, because $\alpha_1\circ d_e$ is continuous and positive definite, $\lim_{k\to\infty} f^{k}(x_0) = x_e$. \qed
	\end{enumerate}

	\begin{rem}
		Note that, because $f$ is a function whose domain and codomain are equal to $M$, $f^k$, where $k \in \mathbb Z_+$, is also a function whose domain and codomain are equal to $M$. Therefore, the existence and uniqueness of the sequence $\{f^k(x_0)\}_{k\in\mathbb Z_+}$ is guaranteed for any $x_0 \in M$.
	\end{rem}
	
	
	\textit{Proof of Theorem \ref{thm:mpc}.}~Given $x_k$, denote by $x_{k+i|k}^*$ the predicted state at time $k+i$ by (\ref{equ:mpcopt}) at time $k$. 
	
	Assume $x_k \in \mathcal D_N$ so that there exists an optimal sequence $\{u^*_{k|k},\dots,u^*_{k+N-1|k}\}$ solving (\ref{equ:mpcopt}). Together, $x_{k+1} = f(x_k,u^*_{k|k})$ and $x^*_{k+N|k} \in \mathcal X_T$ imply that $\{u^*_{k+1|k},\dots,u^*_{k+N-1|k},\kappa(x^*_{k+N|k})\}$ is feasible for (\ref{equ:mpcopt}) at time $k+1$, which implies that $x_{k+1} \in \mathcal D_N$. Since $x_0 \in \mathcal D_N$, then by induction, $x_k \in \mathcal D_N$ for all $k \in \mathbb Z_+$, proving (i).
	
	For $x_k \in \mathcal X_T$, let $x_{k+i+1|k} = f(x_{k+i|k},\kappa(x_{k+i|k}))$ for $i \in \mathbb Z_N$ and $x_{k|k} = x_k$. Let $J_\kappa: \mathcal X_T \rightarrow \mathbb R$ be a function such that $J_\kappa(x_k) = \mathcal V_N(x_k,\{\kappa(x_{k+i|k})\}_{k\in \mathbb Z_N})$. The function $J_\kappa$ is positive definite on $\mathcal X_T$ and continuous, because $\kappa$ and $\mathcal V_N$ are continuous.
	
	Let $J: \mathcal D_N \rightarrow \mathbb R$ be a function such that $J(x_k) = \mathcal V_N^*(x_k)$. Let $J^* = J|_{\mathcal X_T}$ be the restriction of $J$ to the domain $\mathcal X_T \subset \mathcal D_N$. Let $\alpha_1 = \gamma+N\alpha$. 
	To prove that $J^*$ is a Lyapunov function we use the following result.
	
	\begin{lem}
		Let $J_\kappa: M \rightarrow \mathbb R$ be a function that is continuous and positive definite on a compact neighborhood $V \subset M$ of $x_e$. Then there exists a class-$\mathcal K$ function $\alpha_2$ such that $J_\kappa(x) \leq \alpha_2(d_e(x))$ for all $x \in V$.
	\end{lem}
	
	\textit{Proof.} Let $J_m(\delta) = \max_{x \in V}\{J_\kappa(x): d_e(x) \leq \delta\}$. By construction, $J_m$ is a continuous and non-decreasing function with domain $[0,a_{N_a}]$, where $0 < a_{N_a} < \infty$, and satisfies $J_m(0) = 0$. This implies that the domain can be decomposed into finitely many closed intervals $[0,a_{N_a}] = \cup_{i \in \mathbb Z_{N_a}} [a_i,a_{i+1}]$, where $J_m(\delta)$ is strictly increasing on $[a_i,a_{i+1}]$ when $i$ is even and $J_m(\delta)$ is constant on $[a_i,a_{i+1}]$ when $i$ is odd. 
	
	Choose a small, positive scalar value $\sigma > 0$. Let $\alpha_2(0) = 0$, $\alpha_2(\delta) = \alpha_2(a_i)-J_m(a_i)+J_m(\delta)$ for $\delta \in [a_i,a_{i+1}]$ when $i$ is even, and $\alpha_2(\delta) = \alpha_2(a_i)+\sigma\frac{\delta-a_i}{a_{i+1}-a_i}$ for $\delta \in [a_i,a_{i+1}]$ when $i$ is odd. By construction, $\alpha_2$ is a class-$\mathcal K$ function and $J_m(\delta) \leq \alpha_2(\delta)$ for all $\delta \in [0,a_{N_a}]$. Because $J_\kappa(x) \leq J_m(d_e(x)) \leq \alpha_2(d_e(x))$ for all $x \in V$, the proof is complete.  \qed
	
	According to the lemma, there exists a function $\alpha_2$ such that $\alpha_1(d_e(x)) \leq J^*(x) \leq J_\kappa(x) \leq \alpha_2(d_e(x))$ for all $x \in \mathcal X_T$, where $\alpha_1$ and $\alpha_2$ are class-$\mathcal K$ functions.
	
	

	
	By optimality, $J(x_{k+1})-J(x_k) \leq -F(x_{k+N|k}^*)-L(x_k,u_k)+F(x_{k+N+1|k}^*)+L(x_{k+N|k}^*,\kappa(x_{k+N|k}^*)) \leq -L(x_k,u_k) \leq -\gamma(d(x_k,x_e))$
	for any $x_k \in \mathcal D_N$. 
	This implies that $\gamma(d_e(x_k)) \rightarrow 0$ as $k \rightarrow \infty$, and therefore $x_k \rightarrow x_e$ as $k \rightarrow \infty$ for any $x_k \in \mathcal D_N$. Therefore the equilibrium of $x_e$ is asymptotically stable and its domain of attraction is $\mathcal D_N$, proving (ii).
	\qed

	\bibliographystyle{plain}        
	\bibliography{mpc_manifolds}           

\end{document}